\documentclass[11pt]{amsart}
\usepackage{amsfonts,amsmath,amssymb,amsthm}
\usepackage{enumerate}
\usepackage{hyperref}
\usepackage{interval}
\intervalconfig{soft open fences}

\usepackage{amssymb}

\theoremstyle{plain} 
\newtheorem{theorem}{Theorem}

\newtheorem{lemma}[theorem]{Lemma}

\theoremstyle{definition}
\newtheorem{problem}{Problem}

\theoremstyle{remark}
\newtheorem{remark}{Remark}

\providecommand{\ZZ}{\mathbb{Z}}

\providecommand{\NN}{\mathbb{N}}

\providecommand{\mscr}{\mathcal}

\def\ii#1{^{(#1)}}

\def\E{\mathsf{E}}
\def\P{\mathsf{P}}

\let\tl=\tilde

\let\X=X
\let\x=x

\def\X{X}

\def\J{\mathrm{J}}

\def\FK{\mathsf{FK}}

\def\diam{\operatorname{diam}}

\begin{document}
\author{Noam Berger, Anders Johansson, and Anders \"Oberg}
\date{}

\title[Equality of critical inverse temperatures]{Equality of the critical inverse temperatures for the one- and two-sided Dyson models}
\begin{abstract}
For the FK-Ising random-cluster representation of the Dyson long-range Ising
model, we prove that the critical inverse temperatures on the half-line and on
the full line coincide whenever the interaction exponent satisfies
$1<\alpha<2$.  The proof combines a finite-volume cluster estimate above the
full-line critical point, a sprinkling comparison for random-cluster measures,
and a block renormalization that reduces the half-line problem to the classical
one-dimensional long-range site-bond percolation theorem of Newman and
Schulman.  The borderline case $\alpha=2$ remains open for this argument.
\end{abstract}
\maketitle

\section{Introduction}
The Dyson model is the one-dimensional Ising model with pair interaction
$J(ij)=|i-j|^{-\alpha}$.  For $1<\alpha\le2$ the interaction is summable, but
the model remains genuinely long range: the full-line model has a phase
transition at low temperature, and the usual one-dimensional nearest-neighbour
intuition does not apply.

There are two natural geometries.  One may put the model on the full line
$\ZZ$, or on the half-line $\NN$.  In the FK-Ising random-cluster
representation, deleting the negative half-line can only make long connections
harder.  Consequently the half-line critical inverse temperature is at least
as large as the full-line one.  The purpose of this note is to prove the
opposite inequality for $1<\alpha<2$, and hence equality of the two critical
points.

One motivation comes from the transfer-operator approach to one-sided Dyson
potentials.  Johansson, \"Oberg and Pollicott~\cite{jop2025} proved existence
of continuous eigenfunctions in the subcritical regime $\beta<\beta_c^\ZZ$
(when $\alpha>3/2$).  Since a priori one only
knows $\beta_c^\ZZ\le \beta_c^\NN$, the equality proved here shows, for
$1<\alpha<2$, that this full-line threshold is the natural one-sided threshold
as well.  Previously, Johansson, \"Oberg and Pollicott~\cite{jop2019} proved
the quantitative comparison $\beta_c^\NN\le 8\beta_c^\ZZ$.

The proof is based on a renormalization idea 
from
~\cite{berger0}.  Fix
 $\gamma>\alpha/2$. Then above
the full-line critical point, large finite intervals contain,
with high probability,
internal
clusters of size
at least $N^\gamma$.
 Domain monotonicity transfers this estimate to
blocks of the half-line.  After adding an independent long-range Bernoulli
sprinkling, the large block clusters dominate a one-dimensional site-bond
long-range percolation model.  The theorem of Newman and Schulman for
$1/r^\alpha$ percolation with $\alpha<2$ then gives percolation on the
half-line once the block scale is chosen large enough.  A final sprinkling
comparison converts this percolation statement back to the FK-Ising model.

The argument uses the strict inequality $2\gamma-\alpha>0$ with
$\gamma<1$, and therefore leaves the borderline case $\alpha=2$ open.

\medskip

\noindent
{\bf Acknowledgement}. We would like to thank Mirmukhsin Makhmudov for reading
and suggesting improvements of an earlier version of the paper.

\section{The critical inverse temperatures are equal when  \texorpdfstring{$1<\alpha<2$}{1<alpha<2}}

\noindent
Throughout we write $\NN=\{0,1,2,\ldots\}$.  For a finite graph $G$, let
$\mscr C(G)=\{C_1,C_2,\dots\}$ be the set of connected components, let
$\hat C(G)$ denote a component of maximum size, and let
$\omega(G)=|\mscr C(G)|$.

We consider graphs whose vertex sets are integer intervals $V\subset\ZZ$;
$V\ii2$ denotes the set of unordered pairs of distinct vertices in $V$.  The
independent Bernoulli graph model $\eta(G;p)$ is parameterized by an edge
probability function $p:V\ii2\to[0,1]$.  For an interaction strength
$\J(ij)\ge0$, the corresponding Dyson edge probabilities are
\begin{equation}\label{eq:edgeprob}
  p_\beta(ij)=1-\exp(-\beta\J(ij)), \qquad i,j\in\ZZ.
\end{equation}
Here
\begin{equation}\label{eq:dyson}
  \J(ij)=|i-j|^{-\alpha}, \qquad i,j\in\ZZ,
\end{equation}
where $\alpha>1$ and $\beta\ge0$.  We write $p_\beta\vert_V$ for the
restriction of $p_\beta$ to $V\ii2$.

For a finite interval $V\subset\ZZ$, the FK-Ising random-cluster measure
(see \cite{grimmett}) with free boundary condition is
\[
  \nu_\beta^V(G)=\FK_2(G;p_\beta\vert_V)
  \propto 2^{\omega(G)}\prod_{ij\in G}\bigl(e^{\beta\J(ij)}-1\bigr).
\]
For $t\ge0$, we denote by $\eta_t^V$ the independent Bernoulli graph measure on
$V$ with edge probabilities $1-e^{-t\J(ij)}$.  For $V=\NN$ and $V=\ZZ$, the
free-boundary measures $\nu_\beta^V$ are obtained as weak limits along finite
interval exhaustions.  We call $\nu_\beta^\ZZ$ the two-sided model and
$\nu_\beta^\NN$ the one-sided model.

Let $\mathsf{Perc}$ be the event that the graph contains an infinite connected
component.  For $V=\NN$ or $\ZZ$, define
\[
  \beta_c^V=\sup\{\beta\ge0: \nu_\beta^V(\mathsf{Perc})=0\}.
\]
This zero/non-zero threshold is the only property of the critical point used
below.

\begin{theorem}\label{thm:equal}
  We have $\beta_c^\NN=\beta_c^\ZZ$ when $1<\alpha <2$.
\end{theorem}

\begin{problem}\label{prob:alpha2}
  Prove, or disprove, that the same equality $\beta_c^\NN=\beta_c^\ZZ$
  holds in the borderline Dyson case $\alpha=2$.  We expect this to be true;
  the case is important and is not reached by the argument below, which uses
  $2\gamma-\alpha>0$ with $\gamma<1$.
\end{problem}

\subsection{Proof of Theorem~\ref{thm:equal}}

We use stochastic domination for graph measures ordered by inclusion.  Thus
$\nu'\prec\nu$ means that there is a coupling $(H,G)$ with marginals
$\nu'$ and $\nu$ such that $H\subset G$ almost surely.  By Strassen's theorem
\cite{strassen1965}, this is equivalent on finite edge sets to the inequality
$\int f\,d\nu'\le \int f\,d\nu$ for every bounded increasing function $f$.  In
infinite volume we use the corresponding cylinder-event formulation and pass to
limits along finite exhaustions.

We shall use two standard domination facts.  First, suppose that $V\subset W$
and condition on all edges in $W\ii2\setminus V\ii2$.  By the domain Markov
property, the conditional law of the edges in $V$ is the FK-Ising measure on
$V$ with the boundary condition induced by the exterior connections.  Every
such boundary condition dominates the free one.  Hence, by positive
association, if $A$ is any positive-probability event determined by edges of
$W\ii2\setminus V\ii2$, then
\begin{equation}\label{eq:stochdom1}
  \nu_\beta^V \prec
  \mathcal L_{\nu_\beta^W}\bigl(G[V] \mid A\bigr),
\end{equation}
where the conditional law on the right is viewed as a law on graphs on $V$.
Indeed, this law is a mixture of the conditional measures just described.  In
particular, taking $A$ to be the whole sample space shows that the
free-boundary measure increases stochastically when the interval is enlarged.
Second, for every interval $V$ we have the sprinkling relation
\begin{equation}\label{eq:sprinkle}
  \nu_\beta^V(G)
  \prec (\nu_\beta^V\otimes\eta_{\delta}^{V})(G\cup H)
  \prec \nu^V_{\beta+2\delta}(\tl G),
\end{equation}
where $(G,H)\sim\nu_\beta^V\otimes\eta_\delta^V$, i.e. $G$ and $H$ are independent edge configurations on $V$, sampled respectively by 
$\nu_\beta^V$ and $\eta_\delta^V$.  We give the finite-volume
proof of~\eqref{eq:sprinkle} in Appendix~\ref{app:sprinkle}; the infinite-volume
statement follows by weak limits.

By stochastic dominance~\eqref{eq:stochdom1}, we have
$\beta_c^\NN\ge\beta_c^\ZZ$. It remains to prove the reverse inequality. Fix
$\beta_0>\beta_c^\ZZ$. Our goal is to prove that 
$\nu_\beta^{\NN}(\mathsf{Perc})=1$.

To this end choose $\delta>0$ so small that
$\beta:=\beta_0-2\delta>\beta_c^\ZZ$, and let
$G\sim\nu_\beta^\NN$ and $H\sim\eta_\delta^\NN$ be independent. By the
sprinkling relation~\eqref{eq:sprinkle}, it is enough to show that percolation
occurs in $G\cup H$ with positive probability. Monotonicity in $\beta$ then
implies $\beta_0\ge\beta_c^\NN$.

We base the proof on the following lemma. For an integer interval $I$ contained
in the vertex set of $G$, let
\begin{equation}\label{eq:CIdef}
  \mscr C_I(G) = \{C\cap I : C\in \mscr C(G)\}
\end{equation}
and let $\hat C_I(G)$ be an element in $\mscr C_I(G)$ of maximum size.  Thus
paths are allowed to leave $I$ if $G$ has a larger vertex set.  When only the
induced graph on $I$ is meant, we write
$\hat C_I^{\mathrm{int}}(G):=\hat C(G[I])$, where as before $G[I]$ is the restriction of the graph $G$ to the vertex set $I$.
\begin{lemma}\label{lem:clustersize}
  Assume $1<\alpha<2$, $\beta>\beta_c^\ZZ$, and $1>\gamma>\alpha/2$. For any
  $\epsilon>0$ there exist arbitrarily large integers $N$ with the following
  property: for any integer interval $I\subset\NN$ of length $N$,
 \begin{equation}\label{eq:clustersize}
    \nu_\beta^I\bigl(\{G: |\hat C_I(G)| \ge N^\gamma\}\bigr) > 1-\epsilon.
  \end{equation}
\end{lemma}
We postpone the proof of Lemma \ref{lem:clustersize} to Section \ref{sec:pflem2}.
The estimate~\eqref{eq:clustersize} has the following immediate consequence,
which is the form used in the proof of the theorem.  By domain
monotonicity~\eqref{eq:stochdom1},
if $J\supset I$, then, for $G\sim\nu_\beta^J$, the same lower bound holds for
the induced graph $G[I]$ even after conditioning on any positive-probability
event $\mathcal A$ determined by edges outside $I\ii2$:
\[
  \nu_\beta^J\bigl(|\hat C_I^{\mathrm{int}}(G)|\ge N^\gamma
       \mid \mathcal A\bigr)>1-\epsilon.
\]
Consequently, for a sequence of disjoint intervals, the indicators of the
corresponding internal good-block events, by which we mean the events
$|\hat C_I^{\mathrm{int}}(G)|\ge N^\gamma$, can be coupled from below by
independent Bernoulli variables of parameter $1-\epsilon$; this is the standard
sequential coupling consequence of Strassen's theorem.  The terminology is
introduced formally in the proof of Lemma~\ref{lem:clustersize} below.

Given a finite or countable family $\mscr S=\{S_j:j\in U\}$ of disjoint
subsets of an integer interval $J$ indexed by a set $U$, let
$H\sim \eta_\delta^J$ and let $\tl H \sim \mscr H({\mscr S})$ denote the
random graph on vertex set $U$
where an edge $ij\in U\ii2$ is added precisely when $H$ contains an edge
connecting a vertex in $S_i$ with a vertex in $S_j$. Then $\tl H$ is a
Bernoulli graph $\eta(\tl p)$ with edge probabilities $\tl p(ij)$,
$i,j\in U$, given by
\begin{align}
  \tl p(ij)
  &= 1 - \exp\left(-\delta\sum_{x\in S_i,\, y \in S_j} |x-y|^{-\alpha}\right) \nonumber\\
  &\ge 1-\exp\left(-\delta\, |S_i||S_j|\, D_{ij}^{-\alpha}\right),
     \label{eq:edgeprobH}
\end{align}
where $D_{ij}$ is an upper bound on $\diam(S_i\cup S_j)$ and
$\diam A = \max \{|k-l| : k,l\in A\}$. In particular, when the union has
finite diameter,
\begin{equation}\label{eq:edgeprobHmin}
  \tl p(ij) \ge q :=
  1-\exp\left(-\delta\, (\min_{j\in U}|S_j|)^2\,
  \diam\left(\bigcup_{j\in U}S_j\right)^{-\alpha}\right).
\end{equation}

Assume the lemma to be true.  Fix $\gamma\in(\alpha/2,1)$.  By the
one-dimensional long-range site-bond percolation theorem of Newman and
Schulman~\cite{newman}, there are numbers $\lambda_0<1$ and $B<\infty$ such
that the independent site-bond model on $\NN$ with site density
$\lambda\ge\lambda_0$ and bond probabilities at least
$1-\exp(-B|i-j|^{-\alpha})$ percolates with positive probability.  The
half-line formulation follows from the usual two-sided statement by the same
comparison argument that yields the factor-eight bound of~\cite{jop2019}; the
constant is immaterial here because the block parameter below can be made
arbitrarily large.  Choose $\epsilon>0$ so small that
$1-\epsilon\ge\lambda_0$, and then choose from Lemma~\ref{lem:clustersize}
an arbitrarily large $N$ for which $\delta 2^{-\alpha}N^{2\gamma-\alpha}
\ge B$.

Consider the partition
\[
  \mscr I=\{I_k=[kN,(k+1)N) : k\in\NN\}
\]
of $\NN$ into blocks of length $N$.  Here and below
$G\sim\nu_\beta^\NN$ and $H\sim\eta_\delta^\NN$ are independent.  Let
$C_k=\hat C_{I_k}^{\mathrm{int}}(G)$ be a largest component of the induced graph
$G[I_k]$, let
$\tl H = \mscr H( \{C_i : i \in \NN\} )$, and put
\[
  \mscr G = \{ k\in\NN : |C_k| \ge N^{\gamma}\}.
\]
For $i\ne j$, the sets $C_i$ and $C_j$ are contained in an interval
of length at most $(|i-j|+1)N\le 2|i-j|N$. Hence, by
\eqref{eq:edgeprobH}, if $i,j\in\mscr G$, then
\begin{equation}\label{eq:xxp}
  \tl p (ij) \geq 1-e^{-\beta_N|i-j|^{-\alpha}},
\end{equation}
where
\begin{equation}\label{eq:edgeprob2}
  \beta_N \ge \delta \,2^{-\alpha} N^{2\gamma-\alpha}
  \longrightarrow \infty
  \qquad\text{as $N\to\infty$.}
\end{equation}

We now compare this coarse graph with a site-bond long-range percolation model.
The conditional version of the lemma stated above implies that $\mscr G$
stochastically dominates an independent Bernoulli subset of $\NN$ with
retention parameter $\lambda=1-\epsilon$.  Conditional on the retained
vertices, the edges of $\tl H$ dominate independent long-range bonds with
probabilities bounded from below by~\eqref{eq:xxp}.  By the choice of
$\epsilon$ and $N$, this dominating site-bond model is supercritical.  Hence
$\tl H$ percolates with positive probability.

Every infinite component of $\tl H$ gives an infinite component of $G\cup H$.
Therefore $G\cup H$ percolates with positive probability, and the sprinkling
relation~\eqref{eq:sprinkle} implies percolation with positive probability under
$\nu_{\beta_0}^\NN$.  Since $\beta_0>\beta_c^\ZZ$ was arbitrary, this proves
$\beta_c^\NN\le\beta_c^\ZZ$, completing the proof of Theorem~\ref{thm:equal}.
\hfill\qed

\begin{remark}\label{rem:qfk}
The argument extends without essential change to the random-cluster model with
any fixed cluster weight $q\ge1$.  More precisely, let
\[
  \nu_{\beta,q}^V(G)\propto
  q^{\omega(G)}\prod_{ij\in G}\bigl(e^{\beta\J(ij)}-1\bigr)
\]
be the free $q$-random-cluster measure and let $\beta_c^{V,q}$ be its
percolation threshold.  Free-boundary domain monotonicity, the
uniqueness-and-density input in Lemma~\ref{lem:clustersize}, and the block
renormalization are unchanged.  The second sprinkling domination becomes
\[
  (\nu_{\beta,q}^V\otimes\eta_\delta^V)(G\cup H)
  \prec \nu_{\beta+q\delta,q}^V.
\]
The corresponding calculation is given at the end of
Appendix~\ref{app:sprinkle}.  Thus the same proof yields
$\beta_c^{\NN,q}=\beta_c^{\ZZ,q}$ for $1<\alpha<2$.  We have stated the main
result for $q=2$ in order to keep the focus on the Dyson Ising model.
\end{remark}

\subsection{Proof of Lemma~\ref{lem:clustersize}}\label{sec:pflem2}

Write $\beta_\star$ and $\gamma_\star$ for the parameters called $\beta$ and
$\gamma$ in the statement of the lemma.  Choose
\[
  \alpha/2<\gamma_\star<\gamma<1
\]
and choose $\beta_1$ with $\beta_c^\ZZ<\beta_1<\beta_\star$.  Next choose
$\delta_*>0$ so small that $\beta_1+\delta_*<\beta_\star$, and set, for
$n\ge2$,
\[
  \delta_n=\delta_*2^{-n},
  \qquad
  \beta_n=\beta_1+2\sum_{m=2}^{n}\delta_m.
\]
Then $\beta_n<\beta_\star$ for every $n$.  The role of the small increments
$\delta_n$ is to make the sprinkling used at the $n$th renormalization scale
independent of the clusters constructed at the previous scale.

We say that an interval $I$ of length $M=|I|$ is \emph{good} in a graph $G$ if
\begin{equation}\label{eq:gooddef}
  |\hat C_I(G)| \ge |I|^\gamma.
\end{equation}
This is an increasing event.  For an integer $L\ge0$ and a half-open integer
interval $I=[a,b)$, let $I\pm L=[a-L,b+L)$.  We shall construct a fixed
integer $L\ge0$ and an increasing sequence $M_1,M_2,\ldots$ such that, for
every interval $I$ with $|I|=M_n$,
\begin{equation}\label{eq:weakL}
  \nu_{\beta_n}^{I\pm L}\bigl(\text{$I$ is good in $G$}\bigr)>1-\epsilon_n,
\end{equation}
where $\epsilon_n<\epsilon$.

This implies the lemma.  Indeed, choose $n$ arbitrarily large and put
$N=M_n+2L$, with $n$ so large that
\[
  M_n^\gamma\ge (M_n+2L)^{\gamma_\star}=N^{\gamma_\star}.
\]
If $V\subset\NN$ is any interval of length $N$, let $I$ be the interval obtained
from $V$ by removing $L$ sites from each end, so that $V=I\pm L$.  By
\eqref{eq:weakL}, under $\nu_{\beta_n}^{V}$ the interval $V$ contains a cluster
with at least $N^{\gamma_\star}$ vertices with probability at least
$1-\epsilon_n$.  Since $\beta_n<\beta_\star$ and the event is increasing,
monotonicity in $\beta$ gives the same lower bound under $\nu_{\beta_\star}^V$.
As the values $N=M_n+2L$ are arbitrarily large, this proves
\eqref{eq:clustersize}.

Choose a large integer $c_0$ and, for $n\ge2$, define
\begin{equation}\label{eq:cndef}
  c_n=\max\{n^{2/(1-\gamma)}, c_0\},
  \qquad d_n=c_n^{\gamma-1}.
\end{equation}
We take $c_0$ large enough that $d_n\le1/2$ and
$c_n^\gamma-1\ge c_n^\gamma/2$ for all $n\ge2$.  Put
\[
  \epsilon_n=(1+3d_n)\epsilon_{n-1},\qquad n\ge2.
\]
Since $d_n=O(n^{-2})$, the product $\prod_{n\ge2}(1+3d_n)$ is finite.  Choose
\[
  \epsilon_1=\epsilon\Big/\prod_{k=2}^{\infty}(1+3d_k),
\]
so that $\epsilon_n<\epsilon$ for every $n$.

We now choose the base scale $M_1$ and the fixed buffer $L$.  For
$G\sim\nu_{\beta_1}^\ZZ$ with $\beta_1>\beta_c^\ZZ$, the infinite cluster is
almost surely unique and has positive density $\theta>0$; this is the standard
uniqueness and density statement for translation-invariant finite-energy
percolation measures \cite{newman1981}.  By ergodicity, for all sufficiently
large $M$, an interval of length $M$ contains at least $(\theta/2)M$ vertices of
the infinite cluster with probability arbitrarily close to one.  Choose $M_1$
so large that
\[
  M_1^\gamma\le(\theta/2)M_1.
\]
After setting recursively $M_n=M_{n-1}c_n$, we further choose $M_1$ large
enough that, for every $n\ge2$,
\begin{equation}\label{eq:connchoice}
  2\exp\left(\log c_n
        -\frac{\delta_n}{2}c_n^{\gamma-\alpha}M_{n-1}^{2\gamma-\alpha}\right)
  \le \epsilon_{n-1}d_n/2.
\end{equation}
The additional requirement \eqref{eq:connchoice} is possible because
$2\gamma-\alpha>0$, the factor $M_{n-1}^{2\gamma-\alpha}$ grows
super-exponentially in $n$, and $\delta_n=\delta_*2^{-n}$ decays only
exponentially.

For a fixed translate of an interval $I$ of length $M_1$, the event that $I$
contains at least $M_1^\gamma$ vertices of the infinite cluster is the
increasing limit, as $R\to\infty$, of the local event that at least
$M_1^\gamma$ vertices of $I$ are connected to each other by paths contained in
$I\pm R$.  Hence we may choose $R$ so that this local event has probability at
least $1-2\epsilon_1/3$ under $\nu_{\beta_1}^\ZZ$.  By weak convergence of free
finite-volume measures to $\nu_{\beta_1}^\ZZ$, choose $L\ge R$ so large that the
same local event has probability at least $1-\epsilon_1$ under
$\nu_{\beta_1}^{I\pm L}$, uniformly in the translate of $I$.  Thus
\eqref{eq:weakL} holds at scale $n=1$.

We prove \eqref{eq:weakL} by induction.  Assume it holds at scale $n-1$ and
let $I=[a,a+M_n)$ be an interval of length $M_n$.  Subdivide $I$ into the
$c_n$ intervals
\[
  J_j=[a+jM_{n-1},a+(j+1)M_{n-1}),
  \qquad j=0,\ldots,c_n-1.
\]
Let $G\sim\nu_{\beta_{n-1}}^{I\pm L}$.  All probabilities and
expectations up to the sprinkling step are with respect to this law; we denote
them by $\P$ and $\E$, respectively.  For each subinterval $J_j$, set
$\hat C_j=\hat C_{J_j}(G)$.  By the induction hypothesis and domain
monotonicity, the law induced by $\nu_{\beta_{n-1}}^{I\pm L}$ on the larger box
$J_j\pm L$ stochastically dominates $\nu_{\beta_{n-1}}^{J_j\pm L}$.  Hence
\[
  \P\bigl(|\hat C_j|\ge M_{n-1}^\gamma\bigr)>1-\epsilon_{n-1}
  \qquad\text{for every }j.
\]
Let
\[
  \mathcal J=\{0\le j<c_n: |\hat C_j|\ge M_{n-1}^\gamma\},
  \qquad K=|\mathcal J|.
\]
Then
\[
  \E(c_n-K)\le \epsilon_{n-1}c_n.
\]
Since $c_n^\gamma=d_nc_n$, Markov's inequality gives
\begin{equation}\label{eq:markovi}
  \P(K\ge c_n^\gamma)
  \ge 1-\frac{\epsilon_{n-1}}{1-d_n}
  \ge 1-\epsilon_{n-1}(1+2d_n),
\end{equation}
where the last inequality follows from the choice of $c_0$.

Now add an independent sprinkling $H_n\sim\eta_{\delta_n}^{I\pm L}$; from this
point until the end of the induction step, probabilities are taken under the
product law $\nu_{\beta_{n-1}}^{I\pm L}\otimes\eta_{\delta_n}^{I\pm L}$.
Conditional on $G$ and on the selected subinterval clusters
$\mscr S=\{\hat C_j:j\in\mathcal J\}$, the graph induced by the edges of $H_n$
between the sets in $\mscr S$ has independent edges.  On the event
$K\ge c_n^\gamma$, every selected set has size at least $M_{n-1}^\gamma$ and
is contained in $I$, whose diameter is at most $c_nM_{n-1}$.  Therefore, for
every pair of selected subintervals, by~\eqref{eq:edgeprobHmin},
\[
  \tl p(ij)\ge q_n
  :=1-\exp\bigl(-\delta_n c_n^{-\alpha}M_{n-1}^{2\gamma-\alpha}\bigr).
\]
Thus, conditional on $G$ and $K$, the auxiliary graph dominates an
Erd\H{o}s--R\'enyi graph $\mathbb G(K,q_n)$.

We use the standard connectivity estimate for Erd\H{o}s--R\'enyi graphs that,
when $Kq_n$ is large, the probability of disconnection is bounded, up to a
factor $2$, by the probability of having an isolated vertex; see, for example,
Bollob\'as~\cite[Chapter~7]{bollobas}.  Hence
\begin{align*}
  \P(\tl H_n \text{ disconnected}\mid K\ge c_n^\gamma)
  &\le 2c_n\exp\left(-\delta_n(c_n^\gamma-1)c_n^{-\alpha}
        M_{n-1}^{2\gamma-\alpha}\right) \\
  &\le 2\exp\left(\log c_n
        -\frac{\delta_n}{2}c_n^{\gamma-\alpha}M_{n-1}^{2\gamma-\alpha}\right)
   \le \epsilon_{n-1}d_n/2,
\end{align*}
where the last inequality is \eqref{eq:connchoice}.

If $K\ge c_n^\gamma$ and $\tl H_n$ is connected, then the corresponding good
subinterval clusters are all connected in $G\cup H_n$.  Hence
$|\hat C_I(G\cup H_n)|\ge K M_{n-1}^\gamma\ge M_n^\gamma$.  Combining this with
\eqref{eq:markovi} and the preceding disconnection bound gives
\[
  (\nu_{\beta_{n-1}}^{I\pm L}\otimes\eta_{\delta_n}^{I\pm L})
  \bigl(\text{$I$ is good in $G\cup H_n$}\bigr)
  >1-\epsilon_{n-1}(1+3d_n)=1-\epsilon_n.
\]
By the sprinkling relation \eqref{eq:sprinkle}, the law of $G\cup H_n$ is
stochastically dominated by $\nu_{\beta_{n-1}+2\delta_n}^{I\pm L}
=\nu_{\beta_n}^{I\pm L}$.  Since goodness is increasing, \eqref{eq:weakL}
holds at scale $n$.  The induction is complete. \hfill\qed

\appendix
\section{Proof of the sprinkling relation}\label{app:sprinkle}

We prove~\eqref{eq:sprinkle} in finite volume; the infinite-volume statement
follows by taking weak limits along finite intervals. The first domination in
\eqref{eq:sprinkle} is immediate from the coupling $(G,G\cup H)$. We prove the
second one.

Let $E=V\ii2$ and, for an edge $e=xy$, write
\[
  x_e(t):=e^{t\J(e)}-1, \qquad r_e:=e^{\delta\J(e)}-1.
\]
In finite volume the FK-Ising measure is
\[
  \nu_\beta^V(\omega)\propto
  2^{k(\omega)}\prod_{e:\omega_e=1}x_e(\beta),
\]
where $k(\omega)$ is the number of connected components of $\omega$.
Let $\lambda$ be the law of $G\cup H$ when
$(G,H)\sim\nu_\beta^V\otimes\eta_\delta^V$.  For a fixed union configuration
$u$, the event $G\cup H=u$ is the disjoint union, over $g\subset u$, of the
events in which $G=g$ and $H\cup g=u$.  Summing their probabilities gives, up
to a normalising constant,
\[
  \lambda(u)
  \propto
  \prod_{e:u_e=1}r_e
  \sum_{g\subset u}2^{k(g)}
     \prod_{e:g_e=1}\frac{x_e(\beta)e^{\delta\J(e)}}{r_e}.
\]
Fix an edge $e=xy$ and a configuration $u$ on $E\setminus\{e\}$. The ratio of
the conditional weights of $u\cup\{e\}$ and $u$ under $\lambda$ is
\[
  R_e(u)=r_e\left(1+\widetilde x_e
       \pi_u(x\leftrightarrow y)
       +\frac{\widetilde x_e}{2}\pi_u(x\not\leftrightarrow y)\right),
  \qquad
  \widetilde x_e:=\frac{x_e(\beta)e^{\delta\J(e)}}{r_e},
\]
where $\pi_u$ is the FK-Ising measure on the subgraph determined by $u$ with
edge weights $\widetilde x_f$. Under $\nu_{\beta+2\delta}^V$, the corresponding
conditional odds are $x_e(\beta+2\delta)$ if $x$ and $y$ are connected in $u$,
and $x_e(\beta+2\delta)/2$ otherwise.

If $x$ and $y$ are connected in $u$, then
\[
  R_e(u)\le r_e+x_e(\beta)e^{\delta\J(e)}
  =e^{(\beta+\delta)\J(e)}-1
  \le e^{(\beta+2\delta)\J(e)}-1.
\]
If $x$ and $y$ are not connected in $u$, then $\pi_u(x\leftrightarrow y)=0$ and
\[
  R_e(u)=r_e+\frac12 x_e(\beta)e^{\delta\J(e)}
  \le \frac12 x_e(\beta+2\delta),
\]
because, with $a=e^{\beta\J(e)}$ and $b=e^{\delta\J(e)}$,
\[
  (ab^2-1)-\{2(b-1)+(a-1)b\}=(b-1)(ab-1)\ge0.
\]
Thus every one-edge conditional probability under $\lambda$ is bounded above by
the corresponding one-edge conditional probability under $\nu_{\beta+2\delta}^V$.
The FK-Ising single-edge conditional probabilities are increasing in the
outside configuration.  Together with the preceding estimates, this gives the
standard single-edge comparison criterion for stochastic domination: for all
configurations $u\le v$ on
$E\setminus\{e\}$,
\[
  \lambda(\omega_e=1\mid \omega_{E\setminus\{e\}}=u)
  \le
  \nu_{\beta+2\delta}^V(\omega_e=1\mid \omega_{E\setminus\{e\}}=v).
\]
Hence $\lambda\prec\nu_{\beta+2\delta}^V$. This proves the second domination
in \eqref{eq:sprinkle}. \hfill\qed

For completeness, the same calculation also proves the $q$-random-cluster
comparison used in Remark~\ref{rem:qfk}.  Replace $2^{k(g)}$ by $q^{k(g)}$.
Then
\[
  R_e(u)=r_e\left(1+\widetilde x_e
       \pi_u(x\leftrightarrow y)
       +\frac{\widetilde x_e}{q}\pi_u(x\not\leftrightarrow y)\right).
\]
If $x$ and $y$ are connected in $u$, then
$R_e(u)\le x_e(\beta+\delta)\le x_e(\beta+q\delta)$.  If they are not
connected, then
\[
  R_e(u)=r_e+\frac1q x_e(\beta)e^{\delta\J(e)}
  \le \frac1q x_e(\beta+q\delta).
\]
Indeed, with $a=e^{\beta\J(e)}$ and $b=e^{\delta\J(e)}$, the last inequality
is equivalent to
\[
  ab^q-1-q(b-1)-(a-1)b\ge0,
\]
which follows from $a\ge1$, $b\ge1$, and
$b^q-qb+q-1\ge0$.  The latter inequality is the tangent-line inequality for
the convex function $b\mapsto b^q$ at $b=1$.


\begin{thebibliography}{999}

\bibitem{berger0}
N.\ Berger, Transience, recurrence and critical behavior for long-range
percolation, \emph{Comm.\ Math.\ Phys.} \textbf{226} (2002), no.\ 3,
531--558.

\bibitem{jop2019}
A.\ Johansson, A.\ \"Oberg and M.\ Pollicott, Phase transitions in long-range
Ising models and an optimal condition for factors of $g$-measures,
\emph{Ergodic Theory Dynam.\ Systems} \textbf{39} (2019), no.\ 5,
1317--1330.

\bibitem{jop2025}
A.\ Johansson, A.\ \"Oberg and M.\ Pollicott, Continuous eigenfunctions of the
transfer operator for Dyson models, \emph{Mathematische Zeitschrift}
\textbf{310} (2025), no.\ 4, article 62, 1--15.

\bibitem{grimmett}
G.\ R.\ Grimmett, \emph{The Random-Cluster Model}, Grundlehren der
mathematischen Wissenschaften, vol.\ 333, Springer, Berlin, 2006.

\bibitem{newman1981}
C.\ M.\ Newman and L.\ S.\ Schulman, Infinite clusters in percolation models,
\emph{J.\ Stat.\ Phys.} \textbf{26} (1981), 613--628.

\bibitem{newman}
C.\ M.\ Newman and L.\ S.\ Schulman, One-dimensional $1/|j-i|^s$ percolation
models: the existence of a transition for $s\leq2$, \emph{Comm.\ Math.\ Phys.}
\textbf{104} (1986), no.\ 4, 547--571.

\bibitem{bollobas}
B.\ Bollob\'as, \emph{Random Graphs}, 2nd ed., Cambridge Studies in Advanced
Mathematics, vol.\ 73, Cambridge University Press, Cambridge, 2001.

\bibitem{strassen1965}
V.\ Strassen, The existence of probability measures with given marginals,
\emph{Ann.\ Math.\ Statist.} \textbf{36} (1965), 423--439.

\end{thebibliography}
\end{document}